\RequirePackage[2023-06-01]{latexrelease}
\documentclass{birkjour}
\usepackage{amssymb}
\usepackage{url}
 \newtheorem{thm}{Theorem}[section]
 \newtheorem{cor}[thm]{Corollary}
 \newtheorem{lem}[thm]{Lemma}
 \newtheorem{prop}[thm]{Proposition}
 \theoremstyle{definition}
 
 \theoremstyle{remark}

 \numberwithin{equation}{section}

\begin{document}

%
%
%
%
%
%
%
%
%

\title[Density Regularity for Cauchy Integrals]
 {Density Regularity from Derivative Growth\\
of Cauchy-Type Integrals}

\author[O. Volkov]{Oleksandr Volkov}

\address{%
Department of Statistics\\
University of California, Berkeley\\
94720 Berkeley\\
USA}

\email{oleksandr\_volkov@berkeley.edu}



\date{January 1, 2004}

\begin{abstract}
We investigate the relationship between the behavior of the
derivative of a Cauchy-type integral near the curve of
integration and the regularity of its density. We consider
Jordan curves and Jordan arcs composed of finitely many
regular \(C^1\) arcs meeting at nondegenerate angles. For a
Cauchy-type integral with continuous density, we establish
an inverse estimate: if the growth of its derivative off the
curve is controlled by a continuous integrable function of
the distance to the curve that is nonincreasing near zero,
then the modulus of continuity of the density at scale \(h\)
is bounded by the integral of this majorant over \((0,h)\).
As a consequence, we prove that power-law decay of the
derivative as the curve is approached forces the density to
be constant. For Jordan curves, combining the inverse and
direct estimates yields a criterion for H\"older continuity
of the density in terms of the growth of the derivative of
the integral. We then apply this criterion to characterize
the H\"older regularity of analytic functions in bounded
Jordan domains whose boundaries satisfy the same geometric
assumptions. The proof relies on uniformly transverse shifts
of local subarcs and a symmetric jump formula.
\end{abstract}

\maketitle

\section{Introduction}
One of the natural problems in the theory of boundary
behavior of functions is to establish a quantitative
relationship between the smoothness of data prescribed on a
curve and the behavior of the associated analytic functions
near that curve. In the direct setting, the smoothness of
the data is used to estimate the behavior of the analytic
function and its derivatives, whereas in the inverse setting
one seeks to estimate the smoothness of the data from
information about the associated analytic function.

For Cauchy-type integrals, in which the density plays the
role of the data on the curve, these questions have been
studied in a number of works. Salaev~\cite{Salaev1976}
obtained a direct Zygmund-type estimate for the singular
Cauchy integral on an arbitrary closed rectifiable Jordan
curve and, for a certain class of such curves, an inverse
estimate generalizing the Bari--Stechkin theorem. The
dependence of the smoothness of Cauchy-type integrals and
their boundary values on the smoothness of the density and
the geometry of the curve was studied by
Dynkin~\cite{Dynkin1979} and Drobek~\cite{Drobek2013}.
Gerus~\cite{Gerus1982} established the sharpness of an
earlier upper estimate for the moduli of smoothness of the
boundary values of Cauchy-type integrals; in subsequent
works~\cite{Gerus1996,Gerus1998}, he studied moduli of
continuity of the integral and conditions for the existence
of derivatives of its continuous extension.
Plaksa~\cite{Plaksa2025} established a necessary and
sufficient condition for the continuous extension of the
real part of a Cauchy-type integral to an Ahlfors-regular
boundary. For fractional Cauchy transforms, the dependence
of their growth on the modulus of continuity of the data
was studied by Chyzhykov and
Beregova~\cite{ChyzhykovBeregova2019}. Directly related to
the behavior of derivatives off the curve is the work of
Shiga~\cite{Shiga2013}, in which the growth of the \(n\)-th
derivative of the Cauchy integral is estimated in terms of
the modulus of continuity of the density and the distance
to the curve of integration.

In the inverse smoothness results cited above, the hypotheses are
formulated in terms of the smoothness of the singular
integral or of the boundary values of the Cauchy-type
integral, rather than the behavior of its derivative off
the curve. Shiga's estimates, by contrast, are direct: the
modulus of continuity of the density controls the growth of
the derivatives of the integral. This naturally leads to
the converse problem of estimating the modulus of
continuity of the density in terms of the growth of the
derivative of the integral off the curve of integration.

In this paper, we solve this problem for Cauchy-type
integrals with a continuous density \(\varphi\) on \(L\),
where \(L\) is either a Jordan curve or a Jordan arc
composed of finitely many regular \(C^1\) arcs meeting at
nondegenerate angles. Let \(\Phi\) denote the corresponding
Cauchy-type integral, and set
\(\delta(z)=\operatorname{dist}(z,L)\). The main result
states that, for a continuous function
\[
\Omega\colon(0,1]\longrightarrow(0,\infty)
\]
such that \(\Omega\in L^1(0,1)\) and \(\Omega\) is
nonincreasing in some neighborhood of zero, the estimate
\[
|\Phi'(z)|
\leq
\Omega(\delta(z)),
\qquad
z\in\mathbb{C}\setminus L,
\quad
0<\delta(z)\leq1,
\]
implies
\[
\omega(\varphi,h)
\preccurlyeq
\int\limits_0^h\Omega(\tau)\,d\tau,
\qquad
0<h\leq1,
\]
where \(\omega(\varphi,h)\) denotes the modulus of
continuity of \(\varphi\) on \(L\). Here and below,
\(A\preccurlyeq B\) means that \(A\leq CB\) for some
constant \(C>0\) independent of the variables under
consideration. Thus, the behavior of the derivative of the
integral off the curve provides an explicit quantitative
estimate for the modulus of continuity of its density.

In particular, the main result implies that an estimate of
the form
\(
|\Phi'(z)|\preccurlyeq\delta(z)^\varepsilon
\)
for some \(\varepsilon>0\) can hold only for a constant
density. For Jordan curves, combining the inverse and direct
estimates yields an equivalence between the condition
\(\varphi\in H(\mu)\) on \(L\) and the estimate
\(
|\Phi'(z)|\preccurlyeq\delta(z)^{\mu-1}
\),
where \(0<\mu<1\). Applying this criterion together with the
Cauchy integral formula gives the corresponding
characterization of Hölder continuity for functions \(f\)
analytic in a bounded Jordan domain \(G\) and continuous on
\(\overline G\), provided that \(\partial G\) satisfies the
geometric assumptions stated above. For the unit disk, the
resulting characterization recovers the classical
Hardy--Littlewood criterion~\cite{Hardy1932}.

The paper is organized as follows. In
Section~\ref{sec_preliminaries}, we establish geometric
properties of local subarcs, construct their uniform
transverse shifts, and prove the symmetric jump formula. In
Section~\ref{sec_modulus_estimates}, we obtain local and
global estimates for the modulus of continuity of the
density and prove the rigidity result. In
Section~\ref{sec_holder_criterion}, we combine the direct and
inverse estimates to obtain a criterion for the Hölder
continuity of the density. Finally, in
Section~\ref{sec_applications}, we establish the extension
of Hölder continuity from the boundary to the closure of a
domain and apply the resulting criterion to analytic
functions.

\section{Geometric constructions and the symmetric jump formula}
\label{sec_preliminaries}
Let \(L\) be either a Jordan curve or a Jordan arc composed
of finitely many regular \(C^1\) arcs that meet at angles
\(\alpha_j\pi\) at their junction points, where
\(
0<\alpha_j<2,
\) 
\(
j=1,2,\ldots,k.
\)
We call the junction points of adjacent arcs the nodes of
\(L\), and the points of \(L\) that are neither nodes nor
endpoints interior regular points. Fix an orientation on
\(L\). In a neighborhood of an interior regular point, the
local sides lying to the left and to the right as \(L\) is
traversed in the chosen orientation are called, respectively,
the positive and negative local sides of \(L\).

Denote by \(\delta(z)\) the distance from a point
\(z\in\mathbb C\) to \(L\), that is,
\[
\delta(z)=\min_{t\in L}|z-t|.
\]

For a function \(\varphi\in C(L)\), consider the Cauchy-type
integral
\begin{equation}\label{eq_1}
\Phi(z)=\frac{1}{2\pi i}\int\limits_L
\frac{\varphi(t)\,dt}{t-z},
\qquad z\notin L,
\end{equation}
where the integral is taken with respect to the chosen
orientation of \(L\).

\begin{lem}\label{lem_local_geometry}
Let \(L\) satisfy the assumptions stated above. Then there
exist constants
\(
h_1>0
\), 
\(
C_1\geq1,
\)
depending only on \(L\), with the following properties.

If \(t_1,t_2\in L\) and
\[
0<|t_1-t_2|\leq h_1,
\]
then one can choose a local subarc
\(\gamma=\gamma(t_1,t_2)\subset L\) joining \(t_1\) and
\(t_2\), which is contained either in one of the constituent
arcs of \(L\) or in the union of two adjacent constituent
arcs and therefore contains at most one node of \(L\).
Moreover,
\[
\operatorname{length}(\gamma)
\leq
C_1|t_1-t_2|.
\]
If \(\gamma\) contains a node, say \(a\), then one also has
\[
|t_1-a|+|t_2-a|
\leq
C_1|t_1-t_2|.
\]
\end{lem}

\begin{proof}
We first establish a local length estimate on each regular
\(C^1\) arc. Parametrize such an arc by arc length:
\(
t=\lambda(s).
\)
Since the unit tangent vector \(\lambda'(s)\) is continuous,
for every \(s_0\) there exists a neighborhood in which
\[
\operatorname{Re}
\bigl(
\lambda'(s)\overline{\lambda'(s_0)}
\bigr)
\geq
\frac{1}{2}.
\]
Therefore, for any \(s_1<s_2\) lying in this neighborhood,
we have
\[
\begin{aligned}
|\lambda(s_2)-\lambda(s_1)|
&\geq
\operatorname{Re}\left(
\bigl(\lambda(s_2)-\lambda(s_1)\bigr)
\overline{\lambda'(s_0)}
\right)
\\
&=
\int\limits_{s_1}^{s_2}
\operatorname{Re}\left(
\lambda'(s)\overline{\lambda'(s_0)}
\right)\,ds
\geq
\frac{1}{2}|s_2-s_1|.
\end{aligned}
\]

Since the parametrization is by arc length,
\[
|s_2-s_1|
=
\operatorname{length}
\bigl(\lambda([s_1,s_2])\bigr),
\]
and hence
\[
\operatorname{length}
\bigl(\lambda([s_1,s_2])\bigr)
\leq
2|\lambda(s_2)-\lambda(s_1)|.
\]
By the compactness of each arc and the finiteness of the
number of arcs, the sizes of these neighborhoods and the
constant in the resulting estimate can be chosen uniformly.

Consider now a node \(a\) at which two adjacent arcs meet,
and let \(\theta_a\) denote the smaller angle between their
one-sided tangent rays pointing away from \(a\). If the angle
at \(a\) is \(\alpha_j\pi\), then
\[
\theta_a
=
\pi\min\{\alpha_j,2-\alpha_j\}>0.
\]
Choose a number
\(
0<\vartheta_a<\theta_a.
\)
By the continuity of the one-sided tangent vectors, there
exists a neighborhood of \(a\) such that, for any points
\(t_1,t_2\ne a\) lying on different arcs incident to \(a\),
the angle between the vectors
\[
t_1-a
\qquad\text{and}\qquad
t_2-a
\]
is at least \(\vartheta_a\).

Set
\[
x=|t_1-a|,
\qquad
y=|t_2-a|.
\]
By the law of cosines,
\[
|t_1-t_2|^2
\geq
x^2+y^2-2xy\cos\vartheta_a.
\]
Moreover,
\[
x^2+y^2-2xy\cos\vartheta_a
-\sin^2\frac{\vartheta_a}{2}(x+y)^2
=
\cos^2\frac{\vartheta_a}{2}(x-y)^2
\geq0.
\]
Consequently,
\[
|t_1-t_2|
\geq
\sin\frac{\vartheta_a}{2}(x+y),
\]
and therefore
\[
|t_1-a|+|t_2-a|
\leq
\frac{1}{\sin(\vartheta_a/2)}
|t_1-t_2|.
\]
If one of the points \(t_1,t_2\) coincides with \(a\), this
estimate holds directly. Thus, in a sufficiently small
neighborhood of \(a\),
\[
|t_1-a|+|t_2-a|
\leq
C_a|t_1-t_2|,
\]
where \(C_a>0\) depends only on the local geometry of \(L\)
at \(a\).

If the local subarc \(\gamma\) joining \(t_1\) and \(t_2\)
does not contain a node, then the local estimate established
above on a regular \(C^1\) arc gives
\[
\operatorname{length}(\gamma)
\leq
C|t_1-t_2|.
\]
If, on the other hand, \(\gamma\) passes through a node
\(a\), then
\[
\begin{aligned}
\operatorname{length}(\gamma)
&=
\operatorname{length}\bigl(\gamma(t_1,a)\bigr)
+
\operatorname{length}\bigl(\gamma(a,t_2)\bigr)
\\
&\leq
C\bigl(|t_1-a|+|t_2-a|\bigr)
\leq
CC_a|t_1-t_2|,
\end{aligned}
\]
where the estimate obtained above was used in the last
inequality.

It remains to choose all constants uniformly and to justify
the existence of \(h_1\). Parametrize the entire curve \(L\)
by arc length. If \(L\) is a Jordan arc, this gives a
homeomorphism
\(
\sigma\colon[0,\ell]\to L,
\)
whereas if \(L\) is a Jordan curve, it gives a homeomorphism
from the circle of circumference \(\ell\) onto \(L\).

The inverse parametrization is uniformly continuous.
Moreover, the number of nodes is finite. Therefore, one can
choose \(\eta>0\) sufficiently small so that every parameter
interval of length at most \(\eta\) is contained either in
the parameter interval corresponding to a single arc or in
the union of two adjacent such intervals and contains at
most one node. By decreasing \(\eta\) if necessary, we may
also ensure that the local estimates established above hold
on every such interval.

The uniform continuity of the inverse parametrization
implies that there exists \(h_1>0\) such that the conditions
\(
t_1,t_2\in L,
\) 
\(
0<|t_1-t_2|\leq h_1
\)
imply that the distance between the corresponding parameter
values is less than \(\eta\); for a closed curve, the cyclic
distance between the parameter values is used. Therefore,
\(t_1\) and \(t_2\) can be joined by a local subarc
\(\gamma=\gamma(t_1,t_2)\) that is contained either in a
single arc or in the union of two adjacent arcs and contains
at most one node.

Finally, since there are only finitely many regular
\(C^1\) arcs and nodes, we may take the maximum of all the
constants obtained above and choose a single constant
\(C_1\geq1\) such that
\[
\operatorname{length}(\gamma)
\leq
C_1|t_1-t_2|,
\]
and, if \(\gamma\) contains a node \(a\),
\[
|t_1-a|+|t_2-a|
\leq
C_1|t_1-t_2|.
\]
This completes the proof.
\end{proof}

\begin{lem}[Uniform transverse shifts]\label{lem_2}
Let \(h_1\) and \(C_1\) be the constants from
Lemma~\ref{lem_local_geometry}. Then there exist constants
\(
0<h_2\leq h_1,
\) 
\(
0<c_0\leq1,
\)
depending only on \(L\), with the following property.

Let \(t_1,t_2\in L\) be distinct interior regular points,
and suppose that the local subarc
\(\gamma=\gamma(t_1,t_2)\), chosen according to
Lemma~\ref{lem_local_geometry}, contains no nodes. Then, for
every \(q\) satisfying
\[
|t_1-t_2|\leq q\leq h_2,
\]
there exists a unit vector \(v\), transverse to every tangent
line to \(\gamma\), such that
\begin{equation}\label{eq_transverse_distance}
c_0\tau
\leq
\delta(t\pm\tau v)
\leq
\tau,
\qquad
t\in\gamma,
\quad
0<\tau\leq q.
\end{equation}
After replacing \(v\) by \(-v\), if necessary, we may assume
that, for every \(t\in\gamma\) and \(0<\tau\leq q\), the
point \(t+\tau v\) lies on the positive local side of \(L\),
whereas \(t-\tau v\) lies on its negative local side.

In particular, the shifted arcs
\[
l_q^\pm
=
\{t\pm qv:t\in\gamma\}
\]
are contained in \(\mathbb{C}\setminus L\). For each choice
of sign, \(l_q^\pm\) joins \(t_1\pm qv\) to
\(t_2\pm qv\) and satisfies
\[
\operatorname{length}(l_q^\pm)
=
\operatorname{length}(\gamma)
\leq
C_1|t_1-t_2|
\leq
C_1q.
\]
\end{lem}

\begin{proof}
Set
\[
\beta_*=
\begin{cases}
\displaystyle
\frac{\pi}{2}
\min_{1\leq j\leq k}
\{1,\alpha_j,2-\alpha_j\},
& k\geq1,\\[2mm]
\displaystyle
\frac{\pi}{2},
& k=0.
\end{cases}
\]
Choose numbers \(\alpha_0,\beta_0\) such that
\(
0<\alpha_0<\beta_0,
\) 
\(
\alpha_0+\beta_0<\beta_*.
\)

We first justify the uniform control of chord directions in
a neighborhood of each node. Let \(a\) be a node of \(L\).
Parametrize the two arcs incident to \(a\) by arc length
measured from \(a\):
\[
\lambda_j(s)=a+s v_j+r_j(s),
\qquad
r_j(s)=o(s),
\qquad
s\to0+,
\qquad
j=1,2,
\]
where \(v_1\) and \(v_2\) are the one-sided unit tangent
vectors to these arcs at \(a\), directed away from \(a\).

Suppose that the smaller angle between the corresponding
tangent rays is \(2\beta_a\). Then, for \(s,t\geq0\) with
\(s+t>0\), we have
\[
|t v_2-s v_1|^2
=
t^2+s^2-2st\cos(2\beta_a).
\]
Consequently,
\[
|t v_2-s v_1|^2
-
\sin^2\beta_a\,(s+t)^2
=
\cos^2\beta_a\,(t-s)^2
\geq0,
\]
and hence
\[
|t v_2-s v_1|
\geq
\sin\beta_a\,(s+t).
\]

On the other hand,
\[
\lambda_2(t)-\lambda_1(s)
=
t v_2-s v_1+r_2(t)-r_1(s).
\]
Therefore,
\[
\frac{
\bigl|
\lambda_2(t)-\lambda_1(s)-(t v_2-s v_1)
\bigr|
}{
|t v_2-s v_1|
}
\leq
\frac{|r_2(t)|+|r_1(s)|}
{\sin\beta_a\,(s+t)}.
\]
Since \(r_j(s)=o(s)\), the right-hand side tends uniformly
to zero as \(s+t\to0+\). Thus, the directions of the chords
joining points on different arcs incident to the node
converge uniformly to the directions of the corresponding
chords joining points on the tangent rays.

For chords joining points on the same arc, the analogous
statement follows from the continuity of its tangent
direction. Therefore, by shrinking the neighborhood of each
node if necessary, we may ensure that the tangent and chord
directions within each of the two incident arcs, as well as
the directions of chords joining points on different arcs,
differ by at most \(\alpha_0\) from the corresponding
directions determined by the two tangent rays.

Using this property, the uniform continuity of the tangent
directions on the regular \(C^1\) arcs, and the compactness
of \(L\), we choose a finite family of relatively open local
pieces
\(
\Lambda_\nu\subset L,
\) 
\(
\nu=1,\ldots,N,
\)
with the following properties.

If \(\Lambda_\nu\) contains no node, then there exists a
fixed unit tangent direction \(u_\nu\) such that the acute
angle between the line \(\mathbb{R}u_\nu\) and the direction
of any tangent, as well as that of any nonzero chord of
\(\Lambda_\nu\), does not exceed \(\alpha_0\).

If \(\Lambda_\nu\) contains a node \(a\), then it consists
of portions of the two arcs incident to \(a\). Moreover, the
tangent and chord directions within each of these arcs, as
well as the directions of chords joining points on different
arcs, differ by at most \(\alpha_0\) from the corresponding
directions determined by the two one-sided tangent rays at
\(a\).

For each local piece \(\Lambda_\nu\), choose a smaller
relatively open subarc \(U_\nu\) such that
\(
\overline{U_\nu}\subset\Lambda_\nu,
\)
and such that the family
\(
\{U_\nu:\nu=1,\ldots,N\}
\)
covers \(L\). Such a finite cover with the stated inclusions
exists by the compactness of \(L\).

Since \(\overline{U_\nu}\) and
\(L\setminus\Lambda_\nu\) are disjoint compact sets, we have
\[
\operatorname{dist}
\bigl(\overline{U_\nu},L\setminus\Lambda_\nu\bigr)>0.
\]
Since the chosen cover is finite, the constant
\[
d_0:=
\min_{1\leq\nu\leq N}
\operatorname{dist}
\bigl(\overline{U_\nu},L\setminus\Lambda_\nu\bigr)
\]
is positive.

Let \(\rho_0>0\) be a Lebesgue number for the cover
\(\{U_\nu\}_{\nu=1}^N\) of \(L\). Choose \(h_2\), with
\(0<h_2\leq h_1\), sufficiently small that
\[
C_1h_2<\rho_0
\qquad\text{and}\qquad
h_2<\frac{d_0}{2}.
\]

Let \(t_1,t_2\), \(\gamma\), and \(q\) now satisfy the
assumptions of the lemma. By
Lemma~\ref{lem_local_geometry},
\[
\operatorname{length}(\gamma)
\leq
C_1|t_1-t_2|
\leq
C_1q
\leq
C_1h_2
<
\rho_0.
\]
Consequently,
\(
\operatorname{diam}(\gamma)<\rho_0,
\)
and the defining property of the Lebesgue number yields an
index \(\nu\) such that
\(
\gamma\subset U_\nu.
\)

We choose a unit vector \(v\) according to the type of the
local piece \(\Lambda_\nu\).

If \(\Lambda_\nu\) contains no node, set
\(
v\perp u_\nu.
\)
Then the acute angle between the line \(\mathbb{R}v\) and
the direction of any tangent or nonzero chord of
\(\Lambda_\nu\) is at least
\[
\frac{\pi}{2}-\alpha_0
\geq
\beta_*-\alpha_0
>
\beta_0.
\]

Suppose now that \(\Lambda_\nu\) contains a node \(a\).
Since \(\gamma\) contains no node, it lies on one of the two
arcs incident to \(a\). Let \(2\beta_a\) be the smaller angle
between the one-sided tangent rays pointing away from \(a\).
If the angle at the node is \(\alpha_j\pi\), then
\[
\beta_a
=
\frac{\pi}{2}
\min\{\alpha_j,2-\alpha_j\}
\geq
\beta_*.
\]
Choose the direction of \(v\) along the bisector of the
smaller sector determined by these tangent rays.

For the two tangent rays, the acute angle between the line
\(\mathbb{R}v\) and the direction of any chord joining two
points on the same ray or points on different rays is at
least \(\beta_a\). Indeed, after rotating the plane, we may
assume that the tangent rays have directions
\(e^{-i\beta_a}\) and \(e^{i\beta_a}\), while \(v\) points
in the direction of the real axis. A chord joining the
points \(se^{-i\beta_a}\) and \(te^{i\beta_a}\) has direction
vector
\[
te^{i\beta_a}-se^{-i\beta_a}
=
(t-s)\cos\beta_a
+
i(t+s)\sin\beta_a.
\]
Hence the acute angle between this chord and the line
\(\mathbb{R}v\) is at least \(\beta_a\). For chords lying
on a single ray, this follows directly from the choice of
the bisector.

Passing from the tangent rays to the incident \(C^1\) arcs
and using the choice of the local piece \(\Lambda_\nu\), we
obtain that the acute angle between the line
\(\mathbb{R}v\) and any nonzero chord with one endpoint in
\(\gamma\) and the other in \(\Lambda_\nu\) is at least
\[
\beta_a-\alpha_0
\geq
\beta_*-\alpha_0
>
\beta_0.
\]
In particular, \(v\) is transverse to every tangent line to
\(\gamma\).

Thus, in both cases, one can choose a unit vector
\(
v=e^{i\theta}
\)
such that the acute angle between the line
\(\mathbb{R}v\) and any nonzero chord joining a point
\(t\in\gamma\) to a point \(t'\in\Lambda_\nu\) is at least
\(\beta_0\).

We now prove the distance estimate. Let \(t\in\gamma\),
\(0<\tau\leq q\), and \(t'\in\Lambda_\nu\). If \(t'\ne t\),
then the angular estimate established above gives
\[
|t\pm\tau v-t'|
\geq
\tau\sin\beta_0.
\]
For \(t'=t\), the same estimate holds, since
\[
|t\pm\tau v-t|=\tau.
\]

If \(t'\notin\Lambda_\nu\), then the inclusion
\(\gamma\subset U_\nu\) and the definition of \(d_0\) imply
\[
|t-t'|\geq d_0.
\]
Since
\[
0<\tau\leq q\leq h_2<\frac{d_0}{2},
\]
we have
\[
|t\pm\tau v-t'|
\geq
|t-t'|-\tau
\geq
d_0-\tau
\geq
\frac{d_0}{2}
\geq
\tau.
\]
Consequently, for every \(t'\in L\),
\[
|t\pm\tau v-t'|
\geq
c_0\tau,
\]
where
\(
c_0=\sin\beta_0\in(0,1].
\)
Taking the minimum over \(t'\in L\), we obtain
\[
\delta(t\pm\tau v)
\geq
c_0\tau.
\]
On the other hand, since \(t\in L\),
\[
\delta(t\pm\tau v)
\leq
|t\pm\tau v-t|
=
\tau.
\]
This proves estimate~\eqref{eq_transverse_distance}.

Since \(v\) is transverse to every tangent line along the
connected subarc \(\gamma\), the vector \(v\) remains on the
same side of the oriented tangent as one moves along
\(\gamma\). After replacing \(v\) by \(-v\), if necessary,
the points \(t+\tau v\) lie on the positive local side of
\(L\), whereas the points \(t-\tau v\) lie on the negative
local side of \(L\).

The lower bound in~\eqref{eq_transverse_distance} implies
that, for \(0<\tau\leq q\),
\(
t\pm\tau v\notin L.
\)
In particular,
\(
l_q^\pm\subset\mathbb{C}\setminus L.
\)
For each choice of sign, \(l_q^\pm\) is a parallel translate
of \(\gamma\); hence it joins \(t_1\pm qv\) to
\(t_2\pm qv\) and has the same length as \(\gamma\).
Therefore, by Lemma~\ref{lem_local_geometry},
\[
\operatorname{length}(l_q^\pm)
=
\operatorname{length}(\gamma)
\leq
C_1|t_1-t_2|
\leq
C_1q.
\qedhere
\]
\end{proof}

The following symmetric jump formula will serve as our next
analytic tool.

\begin{lem}[Symmetric jump formula]
\label{lem_symmetric_jump}
Let \(t_0\in L\) be an interior regular point of \(L\), and
let \(\tau_0\) be the unit tangent vector to \(L\) at \(t_0\),
oriented consistently with the chosen orientation of \(L\).
Let \(v\) be a unit vector transverse to the tangent line to
\(L\) at \(t_0\), that is,
\(
v\notin\mathbb{R}\tau_0.
\)
Choose the sign of \(v\) so that, for all sufficiently small
\(\varepsilon>0\), the point
\[
t_0+\varepsilon v
\]
lies on the positive local side of \(L\), whereas the point
\[
t_0-\varepsilon v
\]
lies on its negative local side.

Then the Cauchy-type integral defined by~\eqref{eq_1} with
density \(\varphi\in C(L)\) satisfies the symmetric jump
formula
\begin{equation}\label{eq_symmetric_jump}
\lim_{\varepsilon\to0+}
\left[
\Phi(t_0+\varepsilon v)
-
\Phi(t_0-\varepsilon v)
\right]
=
\varphi(t_0).
\end{equation}
\end{lem}

\begin{proof}
Parametrize a neighborhood of \(t_0\) in \(L\) by arc
length:
\[
t=\gamma(s),
\qquad
-r\leq s\leq r,
\]
so that
\[
\gamma(0)=t_0,
\qquad
\gamma'(0)=\tau_0,
\qquad
|\gamma'(s)|=1,
\]
and the direction of the parametrization agrees with the
chosen orientation of \(L\). The number \(r>0\) may be
chosen sufficiently small so that the arc
\[
\Gamma=\gamma([-r,r])
\]
contains neither nodes nor endpoints of \(L\) and
\(
\operatorname{dist}
\bigl(t_0,L\setminus\Gamma\bigr)>0.
\)

Set
\[
z_\varepsilon^+=t_0+\varepsilon v,
\qquad
z_\varepsilon^-=t_0-\varepsilon v.
\]
Then
\[
\begin{split}
\Phi(z_\varepsilon^+)-\Phi(z_\varepsilon^-)
&=
\frac{1}{2\pi i}
\int\limits_{\Gamma}
\varphi(t)
\left(
\frac{1}{t-z_\varepsilon^+}
-
\frac{1}{t-z_\varepsilon^-}
\right)\,dt\\
&\quad+
\frac{1}{2\pi i}
\int\limits_{L\setminus\Gamma}
\varphi(t)
\left(
\frac{1}{t-z_\varepsilon^+}
-
\frac{1}{t-z_\varepsilon^-}
\right)\,dt.
\end{split}
\]

Denote the second integral by \(I_\varepsilon\). Since
\(L\setminus\Gamma\) is separated from \(t_0\), for all
sufficiently small \(\varepsilon>0\) we have
\[
\left|
\frac{1}{t-z_\varepsilon^+}
-
\frac{1}{t-z_\varepsilon^-}
\right|
=
\frac{2\varepsilon}
{|t-z_\varepsilon^+|\,|t-z_\varepsilon^-|}
\leq C\varepsilon,
\qquad
t\in L\setminus\Gamma.
\]
Since \(L\) is rectifiable and \(\varphi\) is bounded on
\(L\), it follows that
\[
I_\varepsilon\longrightarrow0,
\qquad
\varepsilon\to0+.
\]

We now consider the integral over the local arc. Since \(v\)
is transverse to the tangent line to \(L\) at \(t_0\), the
real-linear maps
\[
(s,\varepsilon)\longmapsto
s\tau_0\pm\varepsilon v
\]
are nonsingular. Moreover,
\[
\gamma(s)-t_0=s\tau_0+o(s),
\qquad
s\to0.
\]
Therefore, there exist constants
\[
r>0,
\qquad
\varepsilon_0>0,
\qquad
c>0,
\]
such that
\[
|\gamma(s)-t_0\pm\varepsilon v|
\geq
c\sqrt{s^2+\varepsilon^2},
\qquad
|s|\leq r,
\quad
0<\varepsilon\leq\varepsilon_0.
\]
Decreasing the previously chosen \(r\), if necessary, we
may assume that this estimate holds on the entire arc
\(\Gamma\). Consequently,
\[
\left|
\frac{1}{\gamma(s)-z_\varepsilon^+}
-
\frac{1}{\gamma(s)-z_\varepsilon^-}
\right|
\leq
C\frac{\varepsilon}{s^2+\varepsilon^2},
\qquad
|s|\leq r,
\quad
0<\varepsilon\leq\varepsilon_0.
\]

Write
\[
\varphi(\gamma(s))
=
\varphi(t_0)
+
\bigl(\varphi(\gamma(s))-\varphi(t_0)\bigr).
\]
The local integral can then be written as
\[
\varphi(t_0)J_\varepsilon+R_\varepsilon,
\]
where
\[
J_\varepsilon
=
\frac{1}{2\pi i}
\int\limits_{-r}^{r}
\gamma'(s)
\left(
\frac{1}{\gamma(s)-t_0-\varepsilon v}
-
\frac{1}{\gamma(s)-t_0+\varepsilon v}
\right)\,ds
\]
and
\[
\begin{split}
R_\varepsilon
&=
\frac{1}{2\pi i}
\int\limits_{-r}^{r}
\bigl(\varphi(\gamma(s))-\varphi(t_0)\bigr)
\gamma'(s)\\
&\qquad\qquad{}\times
\left(
\frac{1}{\gamma(s)-t_0-\varepsilon v}
-
\frac{1}{\gamma(s)-t_0+\varepsilon v}
\right)\,ds.
\end{split}
\]

The kernel estimate obtained above gives
\[
|R_\varepsilon|
\leq
C
\int\limits_{-r}^{r}
\frac{\varepsilon
|\varphi(\gamma(s))-\varphi(t_0)|}
{s^2+\varepsilon^2}\,ds.
\]
We show that the right-hand side tends to zero. Indeed, for
every \(\eta>0\), the continuity of \(\varphi\) allows us to
choose \(0<\rho<r\) such that
\[
|\varphi(\gamma(s))-\varphi(t_0)|<\eta,
\qquad
|s|<\rho.
\]
Then
\[
\int\limits_{-\rho}^{\rho}
\frac{\varepsilon
|\varphi(\gamma(s))-\varphi(t_0)|}
{s^2+\varepsilon^2}\,ds
\leq
\eta
\int\limits_{-\infty}^{\infty}
\frac{\varepsilon}{s^2+\varepsilon^2}\,ds
=
\pi\eta.
\]
On the remaining part of the interval,
\[
\int\limits_{\rho\leq|s|\leq r}
\frac{\varepsilon
|\varphi(\gamma(s))-\varphi(t_0)|}
{s^2+\varepsilon^2}\,ds
\leq
C\varepsilon
\int\limits_{\rho\leq|s|\leq r}
\frac{ds}{s^2}
\longrightarrow0.
\]
Since \(\eta>0\) is arbitrary, we obtain
\[
R_\varepsilon\longrightarrow0,
\qquad
\varepsilon\to0+.
\]

It remains to compute the limit of \(J_\varepsilon\). Set
\[
\lambda(s)
=
\overline{\tau_0}\bigl(\gamma(s)-t_0\bigr),
\qquad
w=\overline{\tau_0}v.
\]
Then
\[
\lambda(0)=0,
\qquad
\lambda'(0)=1.
\]
Since \(v\) points toward the positive local side of the
oriented curve, we have
\(
\operatorname{Im}w>0.
\)
After the change of variables \(s=\varepsilon u\), we obtain
\[
J_\varepsilon
=
\frac{1}{2\pi i}
\int\limits_{-r/\varepsilon}^{r/\varepsilon}
\lambda'(\varepsilon u)
\left(
\frac{1}{\lambda(\varepsilon u)/\varepsilon-w}
-
\frac{1}{\lambda(\varepsilon u)/\varepsilon+w}
\right)\,du.
\]

Extend the integrand by zero outside the interval
\(
\left[
-\frac{r}{\varepsilon},
\frac{r}{\varepsilon}
\right].
\)
It is then defined on the entire real line. For each fixed
\(u\in\mathbb R\), this interval contains \(u\) for all
sufficiently small \(\varepsilon>0\), and
\[
\frac{\lambda(\varepsilon u)}{\varepsilon}
\longrightarrow u,
\qquad
\lambda'(\varepsilon u)\longrightarrow1.
\]
Moreover, the local kernel estimate obtained above shows
that, for \(0<\varepsilon\leq\varepsilon_0\), the extended
integrand is bounded in absolute value by the integrable
function
\[
\frac{C}{1+u^2},
\]
independently of \(\varepsilon\). Hence, by the dominated
convergence theorem,
\[
\lim_{\varepsilon\to0+}J_\varepsilon
=
\frac{1}{2\pi i}
\int\limits_{-\infty}^{\infty}
\left(
\frac{1}{u-w}
-
\frac{1}{u+w}
\right)\,du.
\]

To evaluate the last integral, observe that
\[
\frac{1}{u-w}-\frac{1}{u+w}
=
\frac{2w}{u^2-w^2}
=
O(|u|^{-2}),
\qquad
|u|\to\infty.
\]
Therefore, the integral along the upper semicircle tends to
zero as its radius tends to infinity. Since
\(\operatorname{Im}w>0\), the only pole in the upper
half-plane is \(u=w\), and its residue is \(1\). By the
residue theorem,
\[
\int\limits_{-\infty}^{\infty}
\left(
\frac{1}{u-w}
-
\frac{1}{u+w}
\right)\,du
=
2\pi i.
\]
Consequently,
\[
J_\varepsilon\longrightarrow1,
\qquad
\varepsilon\to0+.
\]

Combining the estimates obtained above, we conclude that
\[
\lim_{\varepsilon\to0+}
\left[
\Phi(t_0+\varepsilon v)
-
\Phi(t_0-\varepsilon v)
\right]
=
\varphi(t_0).
\]
This completes the proof.
\end{proof}

\section{Estimates for the modulus of continuity and rigidity
of the density}
\label{sec_modulus_estimates}
We now combine the geometric construction provided by
Lemma~\ref{lem_2} with the symmetric jump formula from
Lemma~\ref{lem_symmetric_jump}.

\begin{prop}[Local estimate on a regular arc]
\label{prop_local_estimate}
Let \(h_2\) be the constant from Lemma~\ref{lem_2}, and let
\(
0<h_0\leq h_2.
\)
Suppose that
\[
\Omega\colon(0,h_0]\longrightarrow(0,\infty)
\]
is a nonincreasing function such that
\(
\Omega\in L^1(0,h_0),
\)
that is,
\[
\int\limits_0^{h_0}\Omega(\tau)\,d\tau<\infty.
\]

Let \(\varphi\in C(L)\), and let \(\Phi\) be the corresponding
Cauchy-type integral defined by~\eqref{eq_1}. Suppose that,
for every \(z\in\mathbb{C}\setminus L\) satisfying
\(
0<\delta(z)\leq h_0,
\)
the estimate
\begin{equation}\label{eq_local_derivative_bound}
|\Phi'(z)|
\leq
\Omega(\delta(z))
\end{equation}
holds.

Let \(t_1,t_2\in L\) be distinct interior regular points,
and suppose that the local subarc
\(\gamma=\gamma(t_1,t_2)\), chosen according to
Lemma~\ref{lem_local_geometry}, contains no nodes. Then, for
every \(q\) satisfying
\(
|t_1-t_2|\leq q\leq h_0,
\)
we have
\begin{equation}\label{eq_local_density_estimate}
|\varphi(t_1)-\varphi(t_2)|
\leq
C
\int\limits_0^q\Omega(\tau)\,d\tau,
\end{equation}
where the constant \(C>0\) depends only on \(L\) and is
independent of \(t_1,t_2\), \(q\), \(\varphi\), and
\(\Omega\).
\end{prop}

\begin{proof}
Let \(q\) satisfy
\(
|t_1-t_2|\leq q\leq h_0.
\)
Since \(h_0\leq h_2\), Lemma~\ref{lem_2} applies to the
points \(t_1,t_2\) and the subarc
\(\gamma=\gamma(t_1,t_2)\). Therefore, there exists a unit
vector \(v\), transverse to every tangent line to
\(\gamma\), such that
\begin{equation}\label{eq_prop_transverse_distance}
c_0\tau
\leq
\delta(t\pm\tau v)
\leq
\tau,
\qquad
t\in\gamma,
\quad
0<\tau\leq q.
\end{equation}
After replacing \(v\) by \(-v\), if necessary, we may assume
that the points \(t+\tau v\) lie on the positive local side
of \(L\), whereas the points \(t-\tau v\) lie on the
negative local side of \(L\).

Fix \(0<\varepsilon<q\). Consider the two shifted arcs
\[
l_q^+
=
\{t+qv:t\in\gamma\},
\qquad
l_q^-
=
\{t-qv:t\in\gamma\}.
\]
By Lemma~\ref{lem_2}, they are contained in
\(\mathbb{C}\setminus L\). The arc \(l_q^+\) joins
\(t_1+qv\) to \(t_2+qv\), while \(l_q^-\) joins
\(t_1-qv\) to \(t_2-qv\), and
\[
\operatorname{length}(l_q^+)
=
\operatorname{length}(l_q^-)
=
\operatorname{length}(\gamma)
\leq
C_1|t_1-t_2|
\leq
C_1q.
\]

Join \(t_1+\varepsilon v\) to
\(t_2+\varepsilon v\) by the path consisting of the line
segment
\[
\{t_1+\tau v:\varepsilon\leq\tau\leq q\},
\]
the arc \(l_q^+\), and the line segment
\[
\{t_2+\tau v:\varepsilon\leq\tau\leq q\}.
\]
This entire path is contained in
\(\mathbb{C}\setminus L\). Similarly,
\(t_1-\varepsilon v\) and \(t_2-\varepsilon v\) are joined
by the path consisting of the two line segments
\[
\{t_j-\tau v:\varepsilon\leq\tau\leq q\},
\qquad
j=1,2,
\]
and the arc \(l_q^-\).

By~\eqref{eq_prop_transverse_distance}, every point \(z\) on
these paths satisfies
\[
0<\delta(z)\leq q\leq h_0.
\]
Hence estimate~\eqref{eq_local_derivative_bound} applies
along both paths.

Since \(\Omega\) is nonincreasing, for \(t\in\gamma\) and
\(0<\tau\leq q\), it follows from
\eqref{eq_prop_transverse_distance} that
\[
|\Phi'(t\pm\tau v)|
\leq
\Omega\bigl(\delta(t\pm\tau v)\bigr)
\leq
\Omega(c_0\tau).
\]
Therefore, estimating the integral of \(\Phi'\) along the
path on the positive side, we obtain
\[
\begin{aligned}
|\Phi(t_1+\varepsilon v)-\Phi(t_2+\varepsilon v)|
&\leq
2\int\limits_\varepsilon^q\Omega(c_0\tau)\,d\tau
+\operatorname{length}(l_q^+)\Omega(c_0q)
\\
&\leq
2\int\limits_\varepsilon^q\Omega(c_0\tau)\,d\tau
+C_1q\Omega(c_0q).
\end{aligned}
\]
Similarly, for the path on the negative side,
\[
|\Phi(t_1-\varepsilon v)-\Phi(t_2-\varepsilon v)|
\leq
2\int\limits_\varepsilon^q
\Omega(c_0\tau)\,d\tau
+
C_1q\Omega(c_0q).
\]

Set
\[
D_j(\varepsilon)
=
\Phi(t_j+\varepsilon v)
-
\Phi(t_j-\varepsilon v),
\qquad
j=1,2.
\]
The vector \(v\) is transverse to the tangent lines to \(L\)
at \(t_1\) and \(t_2\), and, by our choice of its sign,
\(t_j+\varepsilon v\) lies on the positive local side of
\(L\), whereas \(t_j-\varepsilon v\) lies on the negative
local side of \(L\), for \(j=1,2\). Therefore, by
Lemma~\ref{lem_symmetric_jump},
\[
D_j(\varepsilon)\longrightarrow\varphi(t_j),
\qquad
\varepsilon\to0+,
\quad
j=1,2.
\]
Together with the estimates obtained above, this gives
\[
\begin{split}
|\varphi(t_1)-\varphi(t_2)|
&=
\lim_{\varepsilon\to0+}
|D_1(\varepsilon)-D_2(\varepsilon)|\leq
4\int\limits_0^q
\Omega(c_0\tau)\,d\tau
+
2C_1q\Omega(c_0q).
\end{split}
\]

After the change of variables \(u=c_0\tau\), we have
\[
\int\limits_0^q\Omega(c_0\tau)\,d\tau
=
\frac{1}{c_0}
\int\limits_0^{c_0q}\Omega(u)\,du
\leq
\frac{1}{c_0}
\int\limits_0^q\Omega(u)\,du.
\]
Moreover, since \(\Omega\) is nonincreasing,
\[
\int\limits_0^{c_0q}\Omega(u)\,du
\geq
c_0q\,\Omega(c_0q),
\]
and hence
\[
q\Omega(c_0q)
\leq
\frac{1}{c_0}
\int\limits_0^q\Omega(u)\,du.
\]
Thus,
\[
|\varphi(t_1)-\varphi(t_2)|
\leq
\frac{4+2C_1}{c_0}
\int\limits_0^q\Omega(\tau)\,d\tau.
\]
Since \(c_0\) and \(C_1\) depend only on \(L\), this proves
estimate~\eqref{eq_local_density_estimate}.
\end{proof}

We now extend the local estimate from
Proposition~\ref{prop_local_estimate} across the nodes of
\(L\) and obtain a global estimate for the modulus of
continuity of the density.

\begin{thm}\label{th_1}
Let
\[
\Omega\colon(0,1]\longrightarrow(0,\infty)
\]
be a continuous function that is nonincreasing in some
neighborhood of zero and satisfies
\(
\Omega\in L^1(0,1).
\)
If, for every \(z\in\mathbb{C}\setminus L\) such that
\(
0<\delta(z)\leq1,
\)
the estimate
\begin{equation}\label{eq_2}
|\Phi'(z)|
\leq
\Omega(\delta(z))
\end{equation}
holds, then, for every \(0<h\leq1\),
\[
\omega(\varphi,h)
\preccurlyeq
\int_0^h\Omega(\tau)\,d\tau,
\]
where
\[
\omega(\varphi,h)
=
\sup\left\{
|\varphi(t_1)-\varphi(t_2)|:
t_1,t_2\in L,\ |t_1-t_2|\leq h
\right\}
\]
is the modulus of continuity of \(\varphi\) on \(L\).
\end{thm}

The constant \(C>0\) implicit in the notation
\(\preccurlyeq\) in the conclusion of
Theorem~\ref{th_1} may depend on \(L\) and on the choice of
\(r_0>0\) such that \(\Omega\) is nonincreasing on
\((0,r_0]\), but is independent of \(h\) and \(\varphi\).

\begin{proof}
Since \(\Omega\) is nonincreasing in some neighborhood of
zero, there exists \(r_0>0\) such that \(\Omega\) is
nonincreasing on \((0,r_0]\). Choose
\(
h_0=\min\{1,h_2,r_0\},
\)
where \(h_2\) is the constant from Lemma~\ref{lem_2}. Then
\(\Omega\) is nonincreasing on \((0,h_0]\), and the condition
\(\Omega\in L^1(0,1)\) implies that
\(
\int\limits_0^{h_0}\Omega(\tau)\,d\tau<\infty.
\)
Moreover, condition~\eqref{eq_2} implies that
\[
|\Phi'(z)|
\leq
\Omega(\delta(z)),
\qquad
0<\delta(z)\leq h_0.
\]
Thus, Proposition~\ref{prop_local_estimate} applies with
this choice of \(h_0\).

We first extend estimate~\eqref{eq_local_density_estimate}
to the endpoints of regular arcs. Let \(x,y\in L\),
\[
0<|x-y|<h_0,
\]
and suppose that the local subarc \(\gamma(x,y)\), chosen
according to Lemma~\ref{lem_local_geometry}, is contained
in the closure of a single regular \(C^1\) arc and contains
no nodes in its interior.

Choose interior regular points
\(x_n,y_n\in\gamma(x,y)\) converging to \(x\) and \(y\),
respectively, so that the corresponding local subarcs
\(
\gamma(x_n,y_n)
\)
are contained in \(\gamma(x,y)\). For every \(q\) such that
\[
|x-y|<q\leq h_0,
\]
we have, for all sufficiently large \(n\),
\[
|x_n-y_n|\leq q.
\]
By Proposition~\ref{prop_local_estimate},
\[
|\varphi(x_n)-\varphi(y_n)|
\leq
C\int\limits_0^q\Omega(\tau)\,d\tau.
\]
Passing to the limit as \(n\to\infty\) and using the
continuity of \(\varphi\), we obtain
\[
|\varphi(x)-\varphi(y)|
\leq
C\int\limits_0^q\Omega(\tau)\,d\tau.
\]
Finally, letting \(q\) decrease to \(|x-y|\), we obtain
\begin{equation}\label{eq_extended_local_estimate}
|\varphi(x)-\varphi(y)|
\leq
C\int\limits_0^{|x-y|}\Omega(\tau)\,d\tau.
\end{equation}
Here we have used the absolute continuity of the function
\(
q\longmapsto
\int\limits_0^q\Omega(\tau)\,d\tau.
\)

Set
\(
h_*=\frac{h_0}{2C_1},
\)
where \(C_1\geq1\) is the constant from
Lemma~\ref{lem_local_geometry}. We first prove the required
estimate for \(0<h\leq h_*\).

Let \(t_1,t_2\in L\) and
\(
0<|t_1-t_2|\leq h.
\)
Since
\[
h\leq h_*\leq h_0\leq h_2\leq h_1,
\]
Lemma~\ref{lem_local_geometry} shows that \(t_1\) and \(t_2\)
can be joined by a local subarc
\[
\gamma=\gamma(t_1,t_2)\subset L,
\]
which is contained either in a single regular arc or in the
union of two adjacent regular arcs and contains at most one
node.

If \(\gamma\) contains no node, then
Proposition~\ref{prop_local_estimate}, together with the
extension~\eqref{eq_extended_local_estimate}, gives
\[
|\varphi(t_1)-\varphi(t_2)|
\leq
C\int\limits_0^{|t_1-t_2|}
\Omega(\tau)\,d\tau
\leq
C\int\limits_0^h\Omega(\tau)\,d\tau.
\]

Suppose now that \(\gamma\) contains a node \(a\). By
Lemma~\ref{lem_local_geometry},
\[
|t_1-a|+|t_2-a|
\leq
C_1|t_1-t_2|
\leq
C_1h.
\]
In particular,
\[
|t_j-a|
\leq
C_1h
\leq
C_1h_*
=
\frac{h_0}{2},
\qquad
j=1,2.
\]
The subarcs joining \(t_1\) to \(a\) and \(a\) to \(t_2\)
are each contained in the closure of a single regular arc.
Therefore, by~\eqref{eq_extended_local_estimate},
\[
\begin{split}
|\varphi(t_1)-\varphi(t_2)|
&\leq
|\varphi(t_1)-\varphi(a)|
+
|\varphi(a)-\varphi(t_2)|\\
&\leq
C\int\limits_0^{|t_1-a|}
\Omega(\tau)\,d\tau
+
C\int\limits_0^{|t_2-a|}
\Omega(\tau)\,d\tau\\
&\leq
2C\int\limits_0^{C_1h}
\Omega(\tau)\,d\tau.
\end{split}
\]
If one of the points \(t_1,t_2\) coincides with \(a\), then
the corresponding difference is zero, and the estimate
above remains valid.

Since \(C_1h\leq h_0\) and \(\Omega\) is nonincreasing on
\((0,h_0]\), we have
\[
\begin{split}
\int\limits_0^{C_1h}\Omega(\tau)\,d\tau
&=
\int\limits_0^h\Omega(\tau)\,d\tau
+
\int\limits_h^{C_1h}\Omega(\tau)\,d\tau\leq
\int\limits_0^h\Omega(\tau)\,d\tau
+
(C_1-1)h\Omega(h).
\end{split}
\]
Moreover, the monotonicity of \(\Omega\) implies that
\[
h\Omega(h)
\leq
\int\limits_0^h\Omega(\tau)\,d\tau.
\]
Consequently,
\[
\int\limits_0^{C_1h}\Omega(\tau)\,d\tau
\leq
C_1\int\limits_0^h\Omega(\tau)\,d\tau.
\]
Thus, when \(\gamma\) contains a node, we also have
\[
|\varphi(t_1)-\varphi(t_2)|
\leq
C'\int\limits_0^h\Omega(\tau)\,d\tau,
\]
where the constant \(C'>0\) is independent of
\(t_1,t_2\) and \(h\).

Therefore,
\begin{equation}\label{eq_small_h_modulus}
\omega(\varphi,h)
\leq
C'\int\limits_0^h\Omega(\tau)\,d\tau,
\qquad
0<h\leq h_*.
\end{equation}

It remains to consider the case \(h_*<h\leq1\). Since \(L\)
has finite length, there exists an integer \(N\), depending
only on \(L\) and \(h_*\), such that any two points
\(t_1,t_2\in L\) can be joined by a subarc partitioned by
points
\[
t_1=\xi_0,\xi_1,\ldots,\xi_m=t_2,
\qquad
m\leq N,
\]
such that
\(
|\xi_j-\xi_{j-1}|\leq h_*,
\) 
\(
j=1,\ldots,m.
\)
Applying~\eqref{eq_small_h_modulus} to each pair of
consecutive points, we obtain
\[
\begin{split}
|\varphi(t_1)-\varphi(t_2)|
&\leq
\sum_{j=1}^m
|\varphi(\xi_j)-\varphi(\xi_{j-1})|\leq
mC'
\int\limits_0^{h_*}\Omega(\tau)\,d\tau\leq
NC'
\int\limits_0^{h_*}\Omega(\tau)\,d\tau.
\end{split}
\]
Therefore,
\[
|\varphi(t_1)-\varphi(t_2)|
\leq
NC'
\int\limits_0^h\Omega(\tau)\,d\tau.
\]
Taking the supremum over all \(t_1,t_2\in L\) satisfying
\(|t_1-t_2|\leq h\), we obtain
\[
\omega(\varphi,h)
\preccurlyeq
\int\limits_0^h\Omega(\tau)\,d\tau,
\qquad
h_*<h\leq1.
\]

Combining this estimate with~\eqref{eq_small_h_modulus}, we
finally obtain
\[
\omega(\varphi,h)
\preccurlyeq
\int\limits_0^h\Omega(\tau)\,d\tau,
\qquad
0<h\leq1.
\]
This completes the proof.
\end{proof}

The following proposition concerns the case in which the
derivative of the Cauchy-type integral tends to zero as the
curve is approached. Since the function
\(x\mapsto x^\varepsilon\) does not satisfy the monotonicity
assumption of Theorem~\ref{th_1}, this result is proved
separately.

\begin{prop}[Rigidity under decay of the derivative]
\label{prop_rigidity}
If, for some \(\varepsilon>0\),
\[
|\Phi'(z)|
\preccurlyeq
\delta(z)^\varepsilon,
\qquad
z\in\mathbb{C}\setminus L,
\quad
0<\delta(z)\leq1,
\]
then \(\varphi\) is constant on \(L\).
\end{prop}

\begin{proof}
The assumption implies that there exists a constant \(K>0\)
such that
\begin{equation}\label{eq_cor_1_derivative_bound}
|\Phi'(z)|
\leq
K\delta(z)^\varepsilon,
\qquad
z\in\mathbb C\setminus L,
\quad
0<\delta(z)\leq1.
\end{equation}

Set
\(
r_*=\min\{1,h_2\},
\)
where \(h_2\) is the constant from Lemma~\ref{lem_2}. We first
establish a local estimate on a single regular arc.

Let \(t_1,t_2\in L\) be distinct interior regular points, and
suppose that the local subarc
\(
\gamma=\gamma(t_1,t_2),
\)
chosen according to Lemma~\ref{lem_local_geometry}, contains no
nodes. Let
\[
|t_1-t_2|\leq q\leq r_*.
\]
We use the same two paths as in the proof of
Proposition~\ref{prop_local_estimate}. By Lemma~\ref{lem_2},
there exists a unit vector \(v\) such that
\[
\delta(t\pm\tau v)\leq\tau,
\qquad
t\in\gamma,
\quad
0<\tau\leq q,
\]
and
\[
\operatorname{length}(l_q^\pm)
=
\operatorname{length}(\gamma)
\leq
C_1q.
\]

Hence, by \eqref{eq_cor_1_derivative_bound}, the contribution of
the two transverse segments on each side of \(L\) is bounded by
\[
2K\int_\eta^q\tau^\varepsilon\,d\tau,
\qquad
0<\eta<q,
\]
whereas the contribution of the shifted arc is bounded by
\(
KC_1q^{1+\varepsilon}.
\)
Estimating \(\Phi'\) along the two paths, applying the symmetric
jump formula from Lemma~\ref{lem_symmetric_jump}, and letting
\(\eta\to0+\), we obtain
\[
\begin{aligned}
|\varphi(t_1)-\varphi(t_2)|
&\leq
4K\int_0^q\tau^\varepsilon\,d\tau
+
2KC_1q^{1+\varepsilon}
\leq
Cq^{1+\varepsilon},
\end{aligned}
\]
where \(C>0\) depends only on \(K\), \(\varepsilon\), and \(L\).
Taking \(q=|t_1-t_2|\), we conclude that
\begin{equation}\label{eq_cor_1_regular_arc}
|\varphi(t_1)-\varphi(t_2)|
\leq
C|t_1-t_2|^{1+\varepsilon}.
\end{equation}

Estimate \eqref{eq_cor_1_regular_arc} extends to the endpoints
of regular arcs. Indeed, let \(x,y\in L\) satisfy
\[
0<|x-y|<r_*,
\]
and suppose that the subarc joining them is contained in the
closure of a single regular arc and contains no nodes in its
interior. For any \(q\) such that
\[
|x-y|<q\leq r_*,
\]
choose interior regular points \(x_n,y_n\) on this subarc such
that
\[
x_n\to x,
\qquad
y_n\to y.
\]
For all sufficiently large \(n\), we have
\(
|x_n-y_n|\leq q,
\)
and hence the estimate proved above gives
\[
|\varphi(x_n)-\varphi(y_n)|
\leq
Cq^{1+\varepsilon}.
\]
Passing to the limit as \(n\to\infty\), using the continuity of
\(\varphi\), and then letting \(q\) decrease to \(|x-y|\), we
obtain
\begin{equation}\label{eq_cor_1_extended_arc}
|\varphi(x)-\varphi(y)|
\leq
C|x-y|^{1+\varepsilon}.
\end{equation}

Set
\[
h_3=\frac{r_*}{2C_1}.
\]
Let \(t_1,t_2\in L\) satisfy
\(
0<|t_1-t_2|\leq h_3.
\)
By Lemma~\ref{lem_local_geometry}, they can be joined by a local
subarc
\(
\gamma=\gamma(t_1,t_2)
\)
containing at most one node.

If \(\gamma\) contains no node, estimate
\eqref{eq_cor_1_extended_arc} applies directly. Suppose that
\(\gamma\) contains a node \(a\). Then
\[
|t_1-a|+|t_2-a|
\leq
C_1|t_1-t_2|,
\]
and
\[
|t_j-a|
\leq
C_1h_3
=
\frac{r_*}{2},
\qquad
j=1,2.
\]
Applying \eqref{eq_cor_1_extended_arc} on the two regular
subarcs adjacent to \(a\), we obtain
\[
\begin{aligned}
|\varphi(t_1)-\varphi(t_2)|
&\leq
|\varphi(t_1)-\varphi(a)|
+
|\varphi(a)-\varphi(t_2)|
\\
&\leq
C|t_1-a|^{1+\varepsilon}
+
C|t_2-a|^{1+\varepsilon}
\\
&\leq
C\bigl(|t_1-a|+|t_2-a|\bigr)^{1+\varepsilon}
\\
&\leq
CC_1^{1+\varepsilon}
|t_1-t_2|^{1+\varepsilon}.
\end{aligned}
\]
Thus, in both cases, there exists a constant \(C_*>0\),
independent of \(t_1\) and \(t_2\), such that
\begin{equation}\label{eq_cor_1_local_all_points}
|\varphi(t_1)-\varphi(t_2)|
\leq
C_*|t_1-t_2|^{1+\varepsilon},
\qquad
|t_1-t_2|\leq h_3.
\end{equation}

Finally, let \(t_1,t_2\in L\) be arbitrary, and choose a subarc
\(\Gamma\subset L\) joining them. Divide \(\Gamma\) into \(n\)
subarcs of equal length with endpoints
\[
t_1=\tau_0,\tau_1,\ldots,\tau_n=t_2.
\]
Then
\[
|\tau_j-\tau_{j-1}|
\leq
\frac{\operatorname{length}(\Gamma)}{n},
\qquad
j=1,\ldots,n.
\]
For all sufficiently large \(n\), estimate
\eqref{eq_cor_1_local_all_points} applies to every consecutive
pair. Therefore,
\[
\begin{aligned}
|\varphi(t_1)-\varphi(t_2)|
&\leq
\sum_{j=1}^n
|\varphi(\tau_j)-\varphi(\tau_{j-1})|
\\
&\leq
C_*n
\left(
\frac{\operatorname{length}(\Gamma)}{n}
\right)^{1+\varepsilon}
\\
&=
C_*
\frac{\operatorname{length}(\Gamma)^{1+\varepsilon}}
     {n^\varepsilon}.
\end{aligned}
\]
Letting \(n\to\infty\), we obtain
\(
\varphi(t_1)=\varphi(t_2).
\)
Since \(t_1,t_2\in L\) were arbitrary, \(\varphi\) is constant
on \(L\).
\end{proof}

\section{A Hölder continuity criterion for Cauchy-type
integrals}
\label{sec_holder_criterion}

We now return to the Hölder setting. To obtain a two-sided
criterion, we complement the inverse estimate from
Theorem~\ref{th_1} with the corresponding direct estimate
for the derivative of a Cauchy-type integral.

\begin{prop}[Direct estimate for the derivative of a
Cauchy-type integral]
\label{prop_direct_derivative_estimate}
Let \(L\) be a Jordan curve satisfying the assumptions stated
above, and let
\(
0<\mu<1.
\)
Suppose that \(\varphi\in C(L)\) belongs to the Hölder class
\(H(\mu)\) on \(L\), that is,
\[
[\varphi]_{\mu,L}
:=
\sup_{\substack{t,\zeta\in L\\ t\ne\zeta}}
\frac{|\varphi(t)-\varphi(\zeta)|}
{|t-\zeta|^\mu}
<\infty.
\]

Then the Cauchy-type integral \(\Phi\) defined
by~\eqref{eq_1} satisfies
\begin{equation}\label{eq_direct_derivative_estimate}
|\Phi'(z)|
\leq
C[\varphi]_{\mu,L}\,
\delta(z)^{\mu-1},
\qquad
z\in\mathbb{C}\setminus L,
\end{equation}
where the constant \(C>0\) depends only on \(L\) and
\(\mu\), but is independent of \(\varphi\) and \(z\).
\end{prop}

\begin{proof}
We first establish the following geometric estimate. There
exists a constant \(C_L>0\), depending only on \(L\), such
that
\begin{equation}\label{eq_prop_direct_length_ball}
\operatorname{length}\bigl(L\cap B(\zeta,r)\bigr)
\leq
C_Lr,
\qquad
\zeta\in L,
\quad
r>0.
\end{equation}

Indeed, set
\(
\ell=\operatorname{length}(L)
\)
and parametrize \(L\) by arc length:
\[
\sigma\colon[0,\ell]\longrightarrow L,
\qquad
\sigma(0)=\sigma(\ell).
\]
For \(s,s_0\in[0,\ell]\), let
\[
d_\ell(s,s_0)
=
\min\{|s-s_0|,\ell-|s-s_0|\}
\]
denote the cyclic distance between the parameters.

Let \(h_1\) and \(C_1\) be the constants from
Lemma~\ref{lem_local_geometry}. Choose
\[
r_L
=
\min\left\{
h_1,\frac{\ell}{4C_1}
\right\}.
\]
Fix \(\zeta=\sigma(s_0)\in L\), and suppose first that
\(0<r\leq r_L\). If
\(
t=\sigma(s)\in L\cap B(\zeta,r),
\)
then
\(
|t-\zeta|<r\leq h_1.
\)
If \(t=\zeta\), then \(d_\ell(s,s_0)=0\), and the required
conclusion is immediate. Suppose therefore that
\(t\ne\zeta\). By Lemma~\ref{lem_local_geometry}, the points
\(\zeta\) and \(t\) can be joined by a local subarc
\(\gamma(\zeta,t)\) such that
\[
\operatorname{length}\bigl(\gamma(\zeta,t)\bigr)
\leq
C_1|t-\zeta|
<
C_1r.
\]
Consequently,
\[
d_\ell(s,s_0)
\leq
\operatorname{length}\bigl(\gamma(\zeta,t)\bigr)
<
C_1r.
\]
Thus,
\[
L\cap B(\zeta,r)
\subset
\sigma\bigl(
\{s\in[0,\ell]:d_\ell(s,s_0)<C_1r\}
\bigr).
\]
Since \(\sigma\) is parametrized by arc length and
\(C_1r\leq\ell/4\), it follows that
\[
\operatorname{length}\bigl(L\cap B(\zeta,r)\bigr)
\leq
2C_1r,
\qquad
0<r\leq r_L.
\]

If \(r\geq r_L\), then
\[
\operatorname{length}\bigl(L\cap B(\zeta,r)\bigr)
\leq
\ell
\leq
\frac{\ell}{r_L}\,r.
\]
Therefore, estimate~\eqref{eq_prop_direct_length_ball} holds
for every \(r>0\) if we set
\[
C_L
=
\max\left\{
2C_1,\frac{\ell}{r_L}
\right\}.
\]

Now let \(z\in\mathbb{C}\setminus L\). Since \(L\) is
compact, there exists a point \(\zeta\in L\) such that
\(
|z-\zeta|=\delta(z).
\)
Set
\(
d=\delta(z).
\)
Differentiating the Cauchy-type integral, we obtain
\[
\Phi'(z)
=
\frac{1}{2\pi i}
\int\limits_L
\frac{\varphi(t)}{(t-z)^2}\,dt.
\]
Since \(L\) is a Jordan curve,
\[
\int\limits_L\frac{dt}{(t-z)^2}=0.
\]
Consequently,
\begin{equation}\label{eq_prop_direct_subtraction}
\Phi'(z)
=
\frac{1}{2\pi i}
\int\limits_L
\frac{\varphi(t)-\varphi(\zeta)}
     {(t-z)^2}\,dt.
\end{equation}

Suppose first that \(0<d\leq1\). Decompose \(L\) into the
sets
\[
E_0=L\cap B(\zeta,2d)
\]
and
\[
E_k
=
L\cap
\bigl(
B(\zeta,2^{k+1}d)
\setminus B(\zeta,2^kd)
\bigr),
\qquad
k\geq1.
\]

For \(t\in E_0\), we have
\[
|t-\zeta|\leq2d,
\qquad
|t-z|\geq d.
\]
Therefore, the definition of the Hölder seminorm and
estimate~\eqref{eq_prop_direct_length_ball} give
\[
\int\limits_{E_0}
\frac{|\varphi(t)-\varphi(\zeta)|}
     {|t-z|^2}\,|dt|
\leq
[\varphi]_{\mu,L}
\frac{(2d)^\mu}{d^2}
\operatorname{length}
\bigl(L\cap B(\zeta,2d)\bigr)
\leq
C[\varphi]_{\mu,L}d^{\mu-1}.
\]

If \(t\in E_k\), where \(k\geq1\), then
\(
|t-\zeta|<2^{k+1}d
\)
and
\[
|t-z|
\geq
|t-\zeta|-|z-\zeta|
\geq
2^kd-d
\geq
2^{k-1}d.
\]
Consequently,
\[
\begin{split}
\int\limits_{E_k}
\frac{|\varphi(t)-\varphi(\zeta)|}
     {|t-z|^2}\,|dt|
&\leq
[\varphi]_{\mu,L}
\frac{(2^{k+1}d)^\mu}
     {(2^{k-1}d)^2}
\operatorname{length}
\bigl(L\cap B(\zeta,2^{k+1}d)\bigr)\\
&\leq
C[\varphi]_{\mu,L}
(2^kd)^{\mu-1}.
\end{split}
\]
Since \(\mu-1<0\), the series
\(
\sum_{k=1}^{\infty}2^{k(\mu-1)}
\)
converges. Summing the estimates obtained above and
using~\eqref{eq_prop_direct_subtraction}, we obtain
\[
|\Phi'(z)|
\leq
C[\varphi]_{\mu,L}d^{\mu-1},
\qquad
0<d\leq1.
\]

It remains to consider the case \(d\geq1\). For every
\(t\in L\),
\[
|\varphi(t)-\varphi(\zeta)|
\leq
[\varphi]_{\mu,L}|t-\zeta|^\mu
\leq
[\varphi]_{\mu,L}
\operatorname{diam}(L)^\mu.
\]
Therefore, by~\eqref{eq_prop_direct_subtraction},
\[
|\Phi'(z)|
\leq
\frac{[\varphi]_{\mu,L}
      \operatorname{diam}(L)^\mu}{2\pi}
\int\limits_L\frac{|dt|}{|t-z|^2}
\leq
C[\varphi]_{\mu,L}d^{-2}.
\]
Since \(d\geq1\) and \(0<\mu<1\), we have
\(
d^{-2}\leq d^{\mu-1}.
\)
Consequently,
\[
|\Phi'(z)|
\leq
C[\varphi]_{\mu,L}d^{\mu-1},
\qquad
d\geq1.
\]

Combining the two cases and returning to the notation
\(d=\delta(z)\), we obtain
\[
|\Phi'(z)|
\leq
C[\varphi]_{\mu,L}\delta(z)^{\mu-1},
\qquad
z\in\mathbb{C}\setminus L.
\]
The constant \(C>0\) depends only on the geometric
characteristics of \(L\) and on \(\mu\). This completes the
proof.
\end{proof}

Combining the direct estimate from
Proposition~\ref{prop_direct_derivative_estimate} with the
inverse estimate from Theorem~\ref{th_1} yields the following
criterion for the Hölder continuity of the density.

\begin{thm}\label{th_2}
Let \(L\) be a Jordan curve satisfying the assumptions stated
above, and let
\(
0<\mu<1.
\)
For \(\varphi\in C(L)\), the following conditions are
equivalent:

\begin{enumerate}
\item The function \(\varphi\) belongs to the Hölder class
\(H(\mu)\) on \(L\), that is,
\[
|\varphi(t_1)-\varphi(t_2)|
\preccurlyeq
|t_1-t_2|^\mu,
\qquad
t_1,t_2\in L.
\]

\item The derivative of the corresponding Cauchy-type
integral \(\Phi\), defined by \eqref{eq_1}, satisfies
\[
|\Phi'(z)|
\preccurlyeq
\delta(z)^{\mu-1},
\qquad
z\in\mathbb{C}\setminus L.
\]
\end{enumerate}
\end{thm}

\begin{proof}
Necessity follows immediately from
Proposition~\ref{prop_direct_derivative_estimate}.

We prove sufficiency. Suppose that
\[
|\Phi'(z)|
\leq
C_0\delta(z)^{\mu-1},
\qquad
z\in\mathbb{C}\setminus L,
\]
for some constant \(C_0>0\). Set
\[
\Omega(x)=C_0x^{\mu-1},
\qquad
0<x\leq1.
\]
Since \(0<\mu<1\), the function \(\Omega\) is positive,
continuous, and nonincreasing on \((0,1]\), and
\[
\int_0^1\Omega(\tau)\,d\tau
=
C_0\int_0^1\tau^{\mu-1}\,d\tau
=
\frac{C_0}{\mu}
<\infty.
\]
Therefore, Theorem~\ref{th_1} applies to \(\Omega\). Hence,
\[
\omega(\varphi,h)
\preccurlyeq
\int_0^h\Omega(\tau)\,d\tau
=
C_0\int_0^h\tau^{\mu-1}\,d\tau
=
\frac{C_0}{\mu}h^\mu,
\qquad
0<h\leq1.
\]
Thus,
\(
\omega(\varphi,h)
\preccurlyeq
h^\mu,
\) 
\(
0<h\leq1.
\)
Consequently,
\[
|\varphi(t_1)-\varphi(t_2)|
\preccurlyeq
|t_1-t_2|^\mu
\]
for all \(t_1,t_2\in L\) such that
\(
|t_1-t_2|\leq1.
\)

If \(|t_1-t_2|>1\), then the continuity of \(\varphi\) and
the compactness of \(L\) imply that
\[
|\varphi(t_1)-\varphi(t_2)|
\leq
2\|\varphi\|_{\infty,L}
\leq
2\|\varphi\|_{\infty,L}|t_1-t_2|^\mu.
\]
Combining the two cases and increasing the implicit constant
in \(\preccurlyeq\), if necessary, we obtain
\[
|\varphi(t_1)-\varphi(t_2)|
\preccurlyeq
|t_1-t_2|^\mu,
\qquad
t_1,t_2\in L.
\]
Therefore, \(\varphi\in H(\mu)\) on \(L\).
\end{proof}

\section{Application to analytic functions}
\label{sec_applications}

To apply the boundary criterion obtained above to analytic
functions, we need the following auxiliary result.

\begin{prop}
[Extension of Hölder continuity from the boundary to the
closure of a domain]
\label{prop_holder_extension}
Let \(G\subset\mathbb{C}\) be a bounded domain,
\(
L=\partial G,
\) 
\(
0<\mu<1,
\)
and let \(f\) be analytic in \(G\) and continuous on
\(\overline G\).

Suppose that the boundary function
\(
\varphi=f|_L
\)
belongs to the Hölder class \(H(\mu)\) on \(L\), that is,
there exists a constant \(M>0\) such that
\begin{equation}\label{eq_prop_holder_boundary}
|\varphi(\xi_1)-\varphi(\xi_2)|
\leq
M|\xi_1-\xi_2|^\mu,
\qquad
\xi_1,\xi_2\in L,
\end{equation}
and that there exists a constant \(A>0\) such that
\begin{equation}\label{eq_prop_holder_derivative}
|f'(z)|
\leq
A\delta(z)^{\mu-1},
\qquad
z\in G,
\end{equation}
where
\(
\delta(z)=\operatorname{dist}(z,L).
\)

Then \(f\) belongs to the Hölder class \(H(\mu)\) on
\(\overline G\). More precisely, there exists a constant
\(C_\mu>0\), depending only on \(\mu\), such that
\begin{equation}\label{eq_prop_holder_closure}
|f(z_1)-f(z_2)|
\leq
C_\mu(A+M)|z_1-z_2|^\mu,
\qquad
z_1,z_2\in\overline G.
\end{equation}
\end{prop}

\begin{proof}
We first establish an auxiliary estimate comparing the value
of \(f\) at an interior point with its value at a nearest
boundary point.

Let \(z\in G\) and set
\(
d=\delta(z).
\)
Since \(G\) is bounded, its boundary \(L\) is compact.
Therefore, there exists a point \(\xi\in L\) such that
\(
|z-\xi|=d.
\)
Set
\(
u=\frac{z-\xi}{|z-\xi|}
\)
and consider the line segment
\(
w(s)=\xi+su,
\) 
\(
0<s\leq d.
\)
Since
\[
|w(s)-z|=d-s<d
\]
and the open disk \(B(z,d)\) is contained in \(G\), we have
\(
w(s)\in G,
\) 
\(
0<s\leq d.
\)

Moreover,
\[
\delta(w(s))\leq|w(s)-\xi|=s.
\]
On the other hand, for every \(\eta\in L\),
\[
|w(s)-\eta|
\geq
|z-\eta|-|z-w(s)|
\geq
d-(d-s)
=
s.
\]
Consequently,
\begin{equation}\label{eq_prop_holder_segment_distance}
\delta(w(s))=s,
\qquad
0<s\leq d.
\end{equation}

For \(0<\rho<d\), integrating \(f'\) along the line segment
from \(w(\rho)\) to \(z\), we obtain
\[
\begin{split}
|f(z)-f(w(\rho))|
&\leq
\int\limits_\rho^d|f'(w(s))|\,ds\leq
A\int\limits_\rho^d
\delta(w(s))^{\mu-1}\,ds=
A\int\limits_\rho^d s^{\mu-1}\,ds.
\end{split}
\]
Letting \(\rho\to0+\) and using the continuity of \(f\) on
\(\overline G\), we obtain
\begin{equation}\label{eq_prop_holder_nearest_boundary}
|f(z)-\varphi(\xi)|
\leq
\frac{A}{\mu}\delta(z)^\mu.
\end{equation}
For \(z\in L\), estimate~\eqref{eq_prop_holder_nearest_boundary}
also holds if we set \(\xi=z\).

Now let \(z_1,z_2\in\overline G\). Set
\(
r=|z_1-z_2|,
\) 
\(
d_j=\delta(z_j),
\) 
\(
j=1,2.
\)
If \(r=0\), the conclusion is immediate, so we assume that
\(r>0\).

Suppose first that
\(
\max\{d_1,d_2\}>r.
\)
Without loss of generality, assume that \(d_1>r\). Then the
entire line segment \([z_1,z_2]\) is contained in the disk
\(B(z_1,d_1)\) and hence in \(G\). Parametrize this segment
by arc length:
\[
z(s)=z_1+su,
\qquad
0\leq s\leq r,
\]
where
\[
u=\frac{z_2-z_1}{|z_2-z_1|}.
\]
Since the distance function to a closed set is
\(1\)-Lipschitz, we have
\[
\delta(z(s))
\geq
\delta(z_1)-|z(s)-z_1|
=
d_1-s.
\]
Since \(\mu-1<0\), it follows that
\[
\delta(z(s))^{\mu-1}
\leq
(d_1-s)^{\mu-1}.
\]
Therefore,
\[
\begin{split}
|f(z_1)-f(z_2)|
&\leq
\int\limits_0^r|f'(z(s))|\,ds\\
&\leq
A\int\limits_0^r(d_1-s)^{\mu-1}\,ds=
\frac{A}{\mu}
\left(
d_1^\mu-(d_1-r)^\mu
\right).
\end{split}
\]
For \(0<\mu<1\), the function \(x^\mu\) is subadditive, and
hence
\[
d_1^\mu-(d_1-r)^\mu
\leq
r^\mu.
\]
Thus,
\begin{equation}\label{eq_prop_holder_deep_case}
|f(z_1)-f(z_2)|
\leq
\frac{A}{\mu}r^\mu.
\end{equation}

It remains to consider the case
\(
d_1\leq r,
\) 
\(
d_2\leq r.
\)
For \(j=1,2\), choose a point \(\xi_j\in L\) such that
\(
|z_j-\xi_j|=d_j.
\)
If \(z_j\in L\), set \(\xi_j=z_j\). Then
\[
|\xi_1-\xi_2|
\leq
|\xi_1-z_1|+|z_1-z_2|+|z_2-\xi_2|
=
d_1+r+d_2
\leq
3r.
\]
Using~\eqref{eq_prop_holder_nearest_boundary} and the Hölder
estimate~\eqref{eq_prop_holder_boundary}, we obtain
\[
\begin{split}
|f(z_1)-f(z_2)|
&\leq
|f(z_1)-\varphi(\xi_1)|
+
|\varphi(\xi_1)-\varphi(\xi_2)|
+
|\varphi(\xi_2)-f(z_2)|\\
&\leq
\frac{A}{\mu}d_1^\mu
+
M|\xi_1-\xi_2|^\mu
+
\frac{A}{\mu}d_2^\mu\leq
\left(
\frac{2A}{\mu}
+
3^\mu M
\right)r^\mu.
\end{split}
\]

Combining this estimate with
\eqref{eq_prop_holder_deep_case}, we obtain
\[
|f(z_1)-f(z_2)|
\leq
C_\mu(A+M)|z_1-z_2|^\mu,
\qquad
z_1,z_2\in\overline G,
\]
where, for example, we may take
\[
C_\mu=\frac{2}{\mu}+3^\mu.
\]
Therefore, \(f\in H(\mu)\) on \(\overline G\). This completes
the proof.
\end{proof}

We now combine the boundary criterion from
Theorem~\ref{th_2} with Proposition~\ref{prop_holder_extension}
and apply the resulting statement to functions analytic in
a Jordan domain.

\begin{cor}\label{cor_2}
Let \(G\) be a bounded Jordan domain whose boundary
\(
L=\partial G
\)
satisfies the assumptions stated above, and let
\(
0<\mu<1.
\)
For a function \(f\) analytic in \(G\) and continuous on
\(\overline G\), the following conditions are equivalent:

\begin{enumerate}
\item The function \(f\) belongs to the Hölder class
\(H(\mu)\) on \(\overline G\), that is,
\[
|f(z_1)-f(z_2)|
\preccurlyeq
|z_1-z_2|^\mu,
\qquad
z_1,z_2\in\overline G.
\]

\item The derivative \(f'\) satisfies
\[
|f'(z)|
\preccurlyeq
\delta(z)^{\mu-1},
\qquad
z\in G,
\]
where
\(
\delta(z)=\operatorname{dist}(z,L).
\)
\end{enumerate}
\end{cor}

\begin{proof}
Set
\(
\varphi=f|_L
\)
and orient the boundary \(L=\partial G\) positively with
respect to \(G\). Consider the corresponding Cauchy-type
integral
\[
\Phi(z)
=
\frac{1}{2\pi i}
\int\limits_L\frac{\varphi(t)}{t-z}\,dt.
\]

Since \(L\) is a rectifiable Jordan curve and \(f\) is
analytic in \(G\) and continuous on \(\overline G\), the
Cauchy--Goursat theorem for rectifiable Jordan
curves~\cite[p.~540]{Walsh1933} gives
\[
\Phi(z)
=
\begin{cases}
f(z), & z\in G,\\
0, & z\in\mathbb{C}\setminus\overline G.
\end{cases}
\]
Consequently,
\begin{equation}\label{eq_cor_2_cauchy_derivative}
\Phi'(z)
=
\begin{cases}
f'(z), & z\in G,\\
0, & z\in\mathbb{C}\setminus\overline G.
\end{cases}
\end{equation}

We prove necessity. Suppose that \(f\in H(\mu)\) on
\(\overline G\). Then the boundary function
\(\varphi=f|_L\) belongs to \(H(\mu)\) on \(L\). By
Theorem~\ref{th_2},
\[
|\Phi'(z)|
\preccurlyeq
\delta(z)^{\mu-1},
\qquad
z\in\mathbb{C}\setminus L.
\]
In particular, for \(z\in G\),
\eqref{eq_cor_2_cauchy_derivative} gives
\[
|f'(z)|
=
|\Phi'(z)|
\preccurlyeq
\delta(z)^{\mu-1}.
\]
Thus, condition~2 holds.

We prove sufficiency. Suppose that there exists a constant
\(A>0\) such that
\[
|f'(z)|
\leq
A\delta(z)^{\mu-1},
\qquad
z\in G.
\]
In view of~\eqref{eq_cor_2_cauchy_derivative}, we have
\[
|\Phi'(z)|
\leq
A\delta(z)^{\mu-1},
\qquad
z\in\mathbb{C}\setminus L,
\]
because \(\Phi'\) vanishes outside \(\overline G\). By
Theorem~\ref{th_2}, the boundary function
\(\varphi=f|_L\) belongs to \(H(\mu)\) on \(L\). Hence,
there exists a constant \(M>0\) such that
\[
|\varphi(\xi_1)-\varphi(\xi_2)|
\leq
M|\xi_1-\xi_2|^\mu,
\qquad
\xi_1,\xi_2\in L.
\]

Proposition~\ref{prop_holder_extension} now applies to
\(f\). Therefore, there exists a constant \(C_\mu>0\),
depending only on \(\mu\), such that
\[
|f(z_1)-f(z_2)|
\leq
C_\mu(A+M)|z_1-z_2|^\mu,
\qquad
z_1,z_2\in\overline G.
\]
Thus, \(f\in H(\mu)\) on \(\overline G\), and condition~1
holds.

This completes the proof.
\end{proof}

For \(G=\mathbb{D}\), Corollary~\ref{cor_2} recovers the
classical Hardy--Littlewood theorem
\cite[p.~411]{Goluzin1969}; thus, it extends that theorem to
bounded Jordan domains whose boundaries satisfy the
geometric assumptions stated above.

\bigskip
\noindent\textbf{Funding} The authors declare that no funds, grants, or other support were
received during the preparation of this manuscript.

\bigskip
\noindent\textbf{Data Availability Statement}
Data sharing is not applicable to this article as no data were created or analyzed in this study.

\bigskip
\noindent\textbf{Declarations}

\bigskip
\noindent\textbf{Conflict of interest}
The authors have no relevant financial or non-financial interests to disclose.

\end{document}